\documentclass[11pt]{amsart}
\usepackage{amsthm, amssymb, latexsym, amsmath, mathtools, color}
\usepackage[a4paper, total={6.5in, 8in}]{geometry}
\usepackage[dvipsnames,svgnames,table]{xcolor}
\usepackage[colorlinks=true,linkcolor=RoyalBlue,urlcolor=RoyalBlue,citecolor=PineGreen]{hyperref}
\usepackage{tikz}
\usepackage{tikz-cd} 
\usetikzlibrary{arrows.meta}
\usepackage{comment, caption}
\usepackage{graphicx,subcaption}

\usepackage{pgf,tikz,pgfplots}
\pgfplotsset{compat=1.14}
\usepackage{pgf,tikz,pgfplots}
\usepackage{mathrsfs}
\usetikzlibrary{arrows}
\usepackage{mathrsfs}
\usetikzlibrary{arrows}
\usepackage{capt-of}
\usetikzlibrary{decorations.pathreplacing}
\usetikzlibrary{cd}
\usetikzlibrary{positioning}
\usetikzlibrary{calc}
\usepackage{algorithmic}
\usepgfplotslibrary{fillbetween}
\usetikzlibrary{intersections}
\usetikzlibrary{patterns}

\usepackage[nameinlink]{cleveref}

\title{More sum-product type counterexamples: products with shifts and $AA+A$}

\author{Oliver Roche-Newton}
\address{Oliver Roche-Newton \\ Institute for Algebra, Johannes Kepler University Linz, Linz, Austria.}
\email{o.rochenewton@gmail.com}

\author{Carl Schildkraut}
\address{Carl Schildkraut \\ Department of Mathematics, Stanford University, Stanford, USA.}
\email{carlsch@stanford.edu}

\author{Audie Warren}
\address{Audie Warren \\ Johann Radon Institute for Computational and Applied Mathematics, Linz, Austria.}
\email{audie.warren@oeaw.ac.at}

\newcommand{\Z}{\mathbb{Z}}
\newcommand{\R}{\mathbb{R}}

\newcommand{\Q}{\mathbb{Q}}

\newcommand{\cO}{\mathcal{O}}

\DeclarePairedDelimiter\abs{\lvert}{\rvert}

\newtheorem{lemma}{Lemma}

\newtheorem{theorem}{Theorem}
\newtheorem{corollary}{Corollary}

\theoremstyle{remark}

\begin{document}

\begin{abstract}
    Adapting the construction disproving the sum-product conjecture over $\mathbb R$ present in \cite{BloomSawinSchildkrautZhelezov}, we show the existence of a constant $c>0$ and arbitrarily large finite sets $A \subseteq \mathbb R$ such that
    $$|AA+A+A| \ll |A|^{2-c}.$$
    As a corollary, all of the sets $A+A$, $AA$, $(A+1)(A+1)$, $A(A+1)$ and $AA+A$ are of size $O(|A|^{2-c})$ for this construction.
\end{abstract}

\maketitle

\section*{Introduction}

The sum-product problem has seen some remarkable progress recently, with the refutation of the sum-product conjecture over the reals by Bloom, Sawin, Schildkraut and Zhelezov \cite{BloomSawinSchildkrautZhelezov}. The proof uses tools from algebraic number theory and takes inspiration from the OpenAI refutation of the Erd\H{o}s unit distance conjecture, which had appeared just a few days earlier -- see the companion paper \cite{alondistances} for more details, as well as the blog post of Bloom \cite{Bloom_blog} for a nice explanation of both of these proofs. Similar ideas were developed further in subsequent work of Pohoata \cite{Pohoata} to refute a conjecture of Elekes concerning the growth of the set $f(A,A)$ when $f$ is non-degenerate in the sense of the Elekes--R\'{o}nyai problem.

In light of these stunning developments, many previously widely believed conjectures in this area are now being seriously questioned. In this paper, we give new results in this direction.

Let us first consider the set $AA+A$. As the set is formed by a combination of addition and multiplication, this so-called ``expander'' is expected to be large for any $A \subset \mathbb R$. 
The best result in this direction is the bound $|AA+A| \geq |A|^{\frac{3}{2}+\frac{3}{170}-o(1)}$, due to Stevens and Warren \cite{StWa22}. A strong conjecture of Balog \cite{Ba11}, stating that $|AA+A| \geq |A|^2$, was refuted in a paper of Roche-Newton, Ruzsa, Shen and Shkredov \cite{RRSS19} via a construction with $|AA+A| \leq |A|^2/(\log\log |A|)^c$. Many experts in this area expected that the bound $|AA+A| \geq |A|^{2-o(1)}$ should still hold, and this belief was alluded to in \cite{RRSS19}. We show that this is false.

\begin{theorem} \label{thm:AA+A}
     There exists an absolute constant $c>0$ such that there are arbitrarily large finite sets $A \subseteq \mathbb R$ satisfying
    $$|AA+A| \leq |A|^{2-c}.$$
\end{theorem}

Another important principle of sum-product theory is that additive shifts disturb multiplicative structure. This principle manifests itself in lower bounds for the maximum of the size of a product set and a shifted product set. The best result in this direction, due to Stevens and Warren \cite{StWa22}, states that, for any $A \subset \mathbb R$,
\[
\max \{|AA|,|(A+1)(A+1)|\} \geq |A|^{1+\frac{11}{38}-o(1)}.
\]
A folklore conjecture in additive combinatorics is that the exponent in the inequality above can be taken to be arbitrarily close to $2$. We show that this is false. 

\begin{theorem} \label{thm:shifts}
    There exists an absolute constant $c>0$ such that there are arbitrarily large finite sets $A \subseteq \mathbb R$ with
    \[
    \max \{|AA|, |(A+1)(A+1)|\} \leq |A|^{2-c}.
    \]
\end{theorem}

Indeed, Theorems \ref{thm:AA+A} and \ref{thm:shifts} follow as a corollary of the following result.

\begin{theorem}\label{thm:main}
    There exists an absolute constant $c>0$ such that there are arbitrarily large finite sets $A \subseteq \mathbb R$ such that
    $$|AA+A+A| \leq |A|^{2-c}.$$
\end{theorem}

Note that both $A+A$ and $AA$ are subsets of a translate of $AA+A+A$, and so Theorem \ref{thm:main} also recovers the main result of \cite{BloomSawinSchildkrautZhelezov}.

\subsection*{Comparing \texorpdfstring{$\R$}{R} with \texorpdfstring{$\Z$}{Z}}
Given a polynomial $f$ in $k$ variables and a subset $T\subset\mathbb C$, we define the \emph{expanding exponent of $f$ over $T$} to be
\[\inf\left\{c:\abs{f(A,\ldots,A)}\geq\abs{A}^c\text{ for all sufficiently large }A\subset T\right\}.\]
There has been much work towards lower-bounding the expanding exponents of various polynomials; for the history of such results, the reader may consult the ``expanders'' section of \cite{BloomWebsite}. 

\Cref{thm:AA+A} can be viewed as stating that the expanding exponent of $f(x,y,z)=xy+z$ over $\R$ is strictly less than $2$. However, a short argument of Shakan \cite{Shakan} gives that the expanding exponent of this same $f$ over $\Z$ is exactly $2$; moreover, this argument shows that if $A$ is a set of $1$-separated positive real numbers, then $\abs{AA+A}\geq\abs{A}^2$. We thus obtain the following corollary.

\begin{corollary}
    The expanding exponent of $f(x,y,z)=xy+z$ over $\Z$ strictly exceeds the expanding exponent of $f$ over $\R$.
\end{corollary}

\noindent As far as we aware, this is the first instance of a result proving that a given polynomial has different expanding exponents over $\mathbb R$ and $\mathbb Z$. In fact, our proof of \Cref{thm:AA+A} gives an expanding exponent strictly less than $2$ over the algebraic integers. 

\section*{Construction}

We begin by recalling the basic objects forming the construction in \cite{BloomSawinSchildkrautZhelezov}. The construction takes place within a number field $K$ with certain properties, whose existence is guaranteed by a theorem of Martinet \cite{Martinet1978}.

\begin{theorem}\label{thm:martinet}
    There exists an absolute constant $C>0$ such that, for infinitely many $d \in \mathbb N$ there exists a totally real number field $K$ with degree $d$ over $\Q$, and such that $\Delta_K \leq C^d$.
\end{theorem}

We recall that in the above, a \textit{totally real} number field $K$ is a number field such that all field embeddings of $K$ into $\mathbb C$ actually map into $\mathbb R$. The number $\Delta_K$ is the \textit{discriminant} of the number field $K$. From now on $K$ will always be a number field obtained by \Cref{thm:martinet}.

Since $K$ has degree $d$ over $\Q$, there are $d$ field embeddings $K \hookrightarrow \R$. Let us denote them by $\sigma_1,...,\sigma_d$. The image of the ring $\cO_K$ of algebraic integers under the Minkowski embedding $\phi\colon \cO_K \rightarrow \mathbb R^d$ given by
$$\alpha \mapsto (\sigma_1(\alpha),...,\sigma_d(\alpha))$$
yields a rank $d$ lattice $\Lambda$ within $\mathbb R^d$. The density of this lattice is controlled by the discriminant of the field -- specifically, the covolume of $\Lambda$ is $\Delta_K^{1/2}$. We will take a collection of algebraic integers $\alpha$ such that the image $\phi(\alpha)$ within $\Lambda$ lies inside a box of side length $2X$:
$$B^+(X) := \{ \alpha \in \cO_K : |\sigma_i(\alpha)| \leq X \text{ for all }1 \leq i \leq d\}.$$

Upper and lower bounds for the size of $B^+(X)$ were proven in \cite{BloomSawinSchildkrautZhelezov}.

\begin{lemma}[{\cite[Lemma 3.3]{BloomSawinSchildkrautZhelezov}}]\label{lem:addbox}
For $K$ a totally real number field of degree $d$, and for any $X \geq 1$, we have
$$\frac{X^d}{\Delta_K^{1/2}} \leq |B^+(X)| \leq (2X+1)^d.$$
\end{lemma}

In addition to the lattice $\Lambda$, we will need a lattice with multiplicative structure. This is obtained from the group of units of $\cO_K$ by taking the logarithms of the absolute values of the embedding. Under the map $\psi\colon \cO_K^\times \rightarrow \R^d$ given by
$$u \mapsto (\log|\sigma_1(u)|,...,\log|\sigma_d(u)|),$$
the group of units $\cO_K^\times$ is mapped into a lattice of rank $d-1$ lying within the hyperplane $x_1 + x_2 + ... +x_d = 0$. We consider this lattice multiplicative in the sense that multiplication of units within $\cO_K^\times$ corresponds to addition within the lattice. We now take a subset of units which lie within a box inside $\R^d$:
$$B^\times(Y) := \left\{ u \in \cO_K^\times : \abs[\big]{\log|\sigma_i(u)|} \leq Y \text{ for all }1 \leq i \leq d\right\}.$$

We combine Lemma 3.1 and Lemma 3.5 from \cite{BloomSawinSchildkrautZhelezov} to give the following upper and lower bound on $|B^\times(Y)|$.

\begin{lemma}\label{lem:multbox}
    For $K$ a totally real number field of degree $d$, and for any $Y \geq 1$ we have
    $$\frac{Y^{d-1}}{d^{1/2} \Delta_K}\leq |B^\times(Y)| \leq 10(5Y+1)^{d-1}.$$
\end{lemma}

Some key observations: 

\begin{itemize}
    \item The sumset $B^+(X) + B^+(X)$ is small, since it is contained within the box $B^+(2X)$.
    \item The product set $B^\times(Y)B^\times(Y)$ is small, as it is contained within $B^\times(2Y)$.
    \item The multiplicatively structured set $B^\times(Y)$ is contained within the additively structured set $B^+(e^Y)$.
    \item If $u \in B^\times(Y)$, then we also have $u^{-1} \in B^\times(Y)$.
\end{itemize}

Our construction itself is the same as in \cite{BloomSawinSchildkrautZhelezov} -- take $Y$ a sufficiently large integer, take $\epsilon >0$ a sufficiently small absolute constant, and then take $X$ a sufficiently large integer such that $e^{3Y} \leq X$, and define
$$P := X + B^+(\epsilon X), \quad G := B^\times(Y).$$
Now set $A := GP$. By the separation of the unit lattice and the assumption that $\epsilon$ is small enough (see \cite[Lemma~3.4]{BloomSawinSchildkrautZhelezov}, as well as the proof of Lemma 4.1), this product set satisfies $|A| = |G||P|$. The product set $GG$ satisfies
$$|GG| \leq |B^\times(2Y)| \leq 10(10Y+1)^{d-1} \leq C_1^d Y^{d-1} \leq C_2^d \Delta_K |G|$$
for some absolute constants $C_1,C_2 >0$. The key new observation in this proof is that
\begin{equation} \label{containment}
AA+A+A \subset GG(PP+G^3P+G^3P).
\end{equation}
Indeed, an arbitrary element of $AA+A+A$ can be written as
\[
g_1p_1g_2p_2+g_3p_3+g_4p_4=g_1g_2(p_1p_2+g_3g_1^{-1}g_2^{-1}p_3+g_4g_1^{-1}g_2^{-1}p_4),
\]
with $g_i \in G$ and $p_i \in P$ for each $1\leq i\leq 4$. Since $B^\times(Y)$ is closed under inversion, \eqref{containment} follows.

Since we have chosen $X$ to be such that $e^{3Y} \leq X$, the product set $G^3P$ is contained within
\[B^+\left(e^Y\right)^3\cdot B^+((1+\epsilon)X)\subset B^+\left((1+\epsilon)e^{3Y}X\right)\subset B^+((1+\epsilon)X^2),\]
and so the sumset $PP + G^3P+G^3P$ is itself contained within $P' := B^+(3(1+\epsilon)^2X^2)$. We then have, using \Cref{lem:addbox},
\begin{equation}\label{AA+A+A-crude}
|AA+A+A| \leq |GGP'| \leq |GG| |P'| \leq (C_2^d \Delta_K |G|) (6(1+\epsilon)^2X^2+1)^d \leq C_3^d\Delta_K |G| X^{2d}
\end{equation}
for some constant $C_3>0$. Aiming to write \eqref{AA+A+A-crude} in terms of $|A|^2$, we use that $\frac{(\epsilon X)^d}{\Delta_K^{1/2}} \leq |P|$, yielding 
\[
|AA+A+A| \leq C_4^d \Delta_K^2|G||P|^2
\]
for some $C_4>0$. Since $K$ was chosen via \Cref{thm:martinet}, we have $\Delta_K \leq C^d$ for some $C$. We conclude that there exists a constant $C_5 > 0$ with
\begin{equation}\label{eq:AA+A+A-fine}
|AA+A+A| \leq \frac{C_5^d}{|G|} |A|^2.
\end{equation}
We now show that this gives a small power saving. We have from \Cref{lem:multbox} that
$$\frac{C_5^d}{|G|} \leq \frac{C_5^d d^{1/2}\Delta_K}{Y^{d-1}} \leq C_6 \left(\frac{C_6}{Y}\right)^{d-1}$$
for some absolute $C_6>0$. We now assume that $Y$ was chosen sufficiently large so that $\frac{C_6}{Y} \leq \frac{1}{2}$, giving
$$|AA+A+A| \leq \frac{2C_6}{2^{d}} |A|^2.$$
But now we use the fact that, since $X$ and $Y$ are fixed, we have $|A| \leq C_7^d$ for some constant $C_7$ which itself depends on $X$ and $Y$. We can then write
$$\frac{1}{2^d} \leq \frac{1}{|A|^c}$$
where $c = \frac{1}{\log(C_7)}$. This concludes the proof of Theorem \ref{thm:main}, since we now have
$$|AA+A+A| \leq 2 C_6 |A|^{2-c}.$$

\section*{Variants}

A natural modification of the argument from the previous section shows that, for any $k \in \mathbb N$, we can obtain a set $A \subset \mathbb R$ such that
\[
|AA+kA| \leq |A|^{2-c}.
\]
Note that in this statement the constant $c$ depends on $k$, since the modification of the proof requires $X$ to be chosen such that $ e^{kY} \leq X$, and the constant $c$ depends itself on $X$.

Since $(A+k)(A+k) \subset AA+2kA+k^2$, we obtain the following corollary which shows that we can have restricted growth with respect to arbitrarily many integer shifts.

\begin{theorem}
    Let $k \in \mathbb N$. Then there exists $c=c(k) >0$ and arbitrarily large finite sets $A \subseteq \R$ such that
    $$\max_{\lambda \in \{0,1,\dots,k\}} |(A+\lambda)(A+\lambda)| \ll |A|^{2-c}.$$
\end{theorem}

There still remain many variants of the sum-product problem for which it appears that the approach in this paper does not help. An intriguing case is that of the set $A(A+A)$, which is something of a twin to $AA+A$. The current best known lower bound for this problem is the estimate $|A(A+A)| \geq |A|^{\frac{3}{2}+\frac{1}{42}-o(1)}$, due to Bloom \cite{Bl25}. Regarding constructions, the best that we are aware of is simply $A=[n]$, which saves a power of a log factor as a consequence of the solution to the Erd\H{o}s multiplication table problem (see Ford \cite{Fo08}). A result of Murphy, Roche-Newton and Shkredov \cite{MRNS15} proves that $|A(A+A+A+A)| \geq |A|^{2-o(1)}$, and so any construction for this problem needs to distinguish between $A+A$ and $A+A+A+A$, which is something that the recent sum-product constructions do not do.

Another problem which has attracted our interest is the problem of how small $\max\{|A+A|,|A^2+A^2|\}$ can be. (Here, we write $A^2$ for the set of squares of elements of $A$.) We were unable to obtain a power saving for this problem, but we suspect that it is possible.

\textbf{Use of AI disclaimer:} The authors used generative AI to develop the ideas for this note; however, everything written within this note is human written.

\section*{Acknowledgments} Oliver Roche-Newton and Audie Warren were partially supported by the Austrian Science Fund (FWF) project PAT2559123. 
Carl Schildkraut is supported by the National Science Foundation Graduate Research Fellowship Program under Grant No.~DGE-2146755.
We thank Thomas Bloom, Jakob Führer, Michalis Kokkinos, Cosmin Pohoata and Misha Rudnev for helpful discussions.

\bibliography{shifts}

@ARTICLE{BloomSawinSchildkrautZhelezov,
       author = {{Bloom}, Thomas F and {Sawin}, Will and {Schildkraut}, Carl and {Zhelezov}, Dmitrii},
        title = "{The sum-product conjecture is false for real numbers}",
      journal = {arXiv e-prints},
         year = {2026},
        month = {May},
          eid = {arXiv:2605.28781},
        pages = {arXiv:2605.28781},
          doi = {10.48550/arXiv.2605.28781},
archivePrefix = {arXiv},
       eprint = {2605.28781},
 primaryClass = {math.NT},
       adsurl = {https://ui.adsabs.harvard.edu/abs/2026arXiv260528781B}
}

@misc{Bloom_blog,
    author       = "Thomas F. Bloom",
    title        = "Sum-product, unit distances, and number fields",
    howpublished = "https://www.erdosproblems.com/forum/thread/blog:6",
    month        = "May 31,",
    year         = "2026",
    url          = "https://www.erdosproblems.com/forum/thread/blog:6"
}

@ARTICLE{Pohoata,
       author={Cosmin Pohoata},
        title={Split primes and the {E}lekes-{R}\'onyai problem}, 
      journal = {arXiv e-prints},
         year = {2026},
        month = {June},
          eid = {arXiv:2606.13619},
        pages = {arXiv:2606.13619},
          doi = {10.48550/arXiv.2606.13619},
archivePrefix = {arXiv},
       eprint = {2606.13619},
 primaryClass = {math.NT},
}

@ARTICLE{alondistances,
       author={Noga Alon and Thomas F. Bloom and W. T. Gowers and Daniel Litt and Will Sawin and Arul Shankar and Jacob Tsimerman and Victor Wang and Melanie Matchett Wood},
        title={Remarks on the disproof of the unit distance conjecture}, 
      journal = {arXiv e-prints},
         year = {2026},
        month = {May},
          eid = {arXiv:2605.20695},
        pages = {arXiv:2605.20695},
          doi = {10.48550/arXiv.2605.20695},
archivePrefix = {arXiv},
       eprint = {2605.20695},
 primaryClass = {math.CO},
}

@article{Martinet1978,
	risfield_0_da = {1978/02/01},
	author = {Martinet, Jacques},
	doi = {10.1007/BF01389902},
	issn = {1432-1297},
	journal = {Inventiones mathematicae},
	number = {1},
	pages = {65–73},
	title = {Tours de corps de classes et estimations de discriminants},
	volume = {44},
	year = {1978}
}

@article {StWa22,
    AUTHOR = {Stevens, Sophie and Warren, Audie},
     TITLE = {On sum sets and convex functions},
   JOURNAL = {Electron. J. Combin.},
  FJOURNAL = {Electronic Journal of Combinatorics},
    VOLUME = {29},
      YEAR = {2022},
    NUMBER = {2},
     PAGES = {Paper No. 2.18, 19},
      ISSN = {1077-8926},
   MRCLASS = {11B30 (05A20 11B13)},
  MRNUMBER = {4418092},
MRREVIEWER = {Boqing\ Xue},
       DOI = {10.37236/10852},
       URL = {https://doi.org/10.37236/10852},
}

@article {RRSS19,
    AUTHOR = {Roche-Newton, Oliver and Ruzsa, Imre Z. and Shen, Chun-Yen and
              Shkredov, Ilya D.},
     TITLE = {On the size of the set {$AA+A$}},
   JOURNAL = {J. Lond. Math. Soc. (2)},
  FJOURNAL = {Journal of the London Mathematical Society. Second Series},
    VOLUME = {99},
      YEAR = {2019},
    NUMBER = {2},
     PAGES = {477--494},
      ISSN = {0024-6107,1469-7750},
   MRCLASS = {52C10 (11B30)},
  MRNUMBER = {3939264},
MRREVIEWER = {Hans\ Parshall},
       DOI = {10.1112/jlms.12177},
       URL = {https://doi.org/10.1112/jlms.12177},
}

@article {Ba11,
    AUTHOR = {Balog, Antal},
     TITLE = {A note on sum-product estimates},
   JOURNAL = {Publ. Math. Debrecen},
  FJOURNAL = {Publicationes Mathematicae Debrecen},
    VOLUME = {79},
      YEAR = {2011},
    NUMBER = {3-4},
     PAGES = {283--289},
      ISSN = {0033-3883,2064-2849},
   MRCLASS = {11B13 (11B75 11P70)},
  MRNUMBER = {2907965},
MRREVIEWER = {Mei\ Chu\ Chang},
       DOI = {10.5486/PMD.2011.5104},
       URL = {https://doi.org/10.5486/PMD.2011.5104},
}

@MISC {Shakan,
    TITLE = {Sum and product estimate over integers, rationals, and reals},
    AUTHOR = {George Shakan},
    HOWPUBLISHED = {MathOverflow},
    NOTE = {URL:https://mathoverflow.net/q/168844 (version: 2015-04-27)},
    EPRINT = {https://mathoverflow.net/q/168844},
    URL = {https://mathoverflow.net/q/168844}
}

@misc {BloomWebsite,
    AUTHOR = {Bloom, Thomas},
    TITLE = {A history of the sum-product problem},
    YEAR = {2026},
    URL = {http://thomasbloom.org/notes/sumproduct.html},
    NOTE = {Available at \url{http://thomasbloom.org/notes/sumproduct.html}, accessed 25 May 2026},
}

@article {Bl25,
  AUTHOR = {Bloom, Thomas},
  TITLE = {Control and its applications in additive combinatorics},
  JOURNAL = {arXiv:2501.09470},
  YEAR = {2025},
}

@article {MRNS15,
    AUTHOR = {Murphy, Brendan and Roche-Newton, Oliver and Shkredov, Ilya},
     TITLE = {Variations on the sum-product problem},
   JOURNAL = {SIAM J. Discrete Math.},
  FJOURNAL = {SIAM Journal on Discrete Mathematics},
    VOLUME = {29},
      YEAR = {2015},
    NUMBER = {1},
     PAGES = {514--540},
      ISSN = {0895-4801,1095-7146},
       DOI = {10.1137/140952004},
       URL = {https://doi.org/10.1137/140952004},
}

@article {Fo08,
    AUTHOR = {Ford, Kevin},
     TITLE = {The distribution of integers with a divisor in a given
              interval},
   JOURNAL = {Ann. of Math. (2)},
  FJOURNAL = {Annals of Mathematics. Second Series},
    VOLUME = {168},
      YEAR = {2008},
    NUMBER = {2},
     PAGES = {367--433},
      ISSN = {0003-486X,1939-8980},
   MRCLASS = {11N25 (11N37)},
  MRNUMBER = {2434882},
MRREVIEWER = {D.\ R.\ Heath-Brown},
       DOI = {10.4007/annals.2008.168.367},
       URL = {https://doi.org/10.4007/annals.2008.168.367},
}
\bibliographystyle{plain}
\end{document}